\documentclass[12pt,reqno]{article}

\usepackage[usenames]{color}
\usepackage{amssymb}
\usepackage{amsmath}
\usepackage{amsthm}
\usepackage{amsfonts}
\usepackage{amscd}
\usepackage{graphicx}

\usepackage[colorlinks=true,
linkcolor=webgreen,
filecolor=webbrown,
citecolor=webgreen]{hyperref}

\definecolor{webgreen}{rgb}{0,.5,0}
\definecolor{webbrown}{rgb}{.6,0,0}

\usepackage{color}
\usepackage{fullpage}
\usepackage{float}

\usepackage{graphics}
\usepackage{latexsym}
\usepackage{epsf}
\usepackage{breakurl}

\newcommand{\seqnum}[1]{\href{https://oeis.org/#1}{\textrm{\underline{#1}}}}
\newcommand{\udiv}{\,{\mid}^{\ast}}

\usepackage{enumerate}
\usepackage{url}

\newcommand{\reg}{\text{reg}}
\newcommand{\znr}{\mathbb{Z}_n^\reg}

\begin{document}

\begin{center}
\epsfxsize=4in
\end{center}

\theoremstyle{plain}
\newtheorem{theorem}{Theorem}
\newtheorem{corollary}[theorem]{Corollary}
\newtheorem{lemma}[theorem]{Lemma}
\newtheorem{proposition}[theorem]{Proposition}

\theoremstyle{definition}
\newtheorem{definition}[theorem]{Definition}
\newtheorem{example}[theorem]{Example}
\newtheorem{conjecture}[theorem]{Conjecture}

\theoremstyle{remark}
\newtheorem{remark}[theorem]{Remark}

\begin{center}
\vskip 1cm{\LARGE\bf On the Number of Regular Integers Modulo $\pmb{n}$ \\
\vskip .1in
and Its Significance for Cryptography
}
\vskip 1cm
\large
Klaus Dohmen and Mandy Lange-Geisler \\
  Fachgruppe Mathematik \\ Hochschule Mittweida
  \\ Technikumplatz 17\\ 09648 Mittweida\\
  Germany \\
\href{mailto: dohmen@hs-mittweida.de}{\tt dohmen@hs-mittweida.de}\\
\href{mailto: mlange1@hs-mittweida.de}{\tt mlange1@hs-mittweida.de}\\
\end{center}

\vskip .2 in
\begin{abstract}
We present four combinatorial proofs of Morgado's formula for the number
$\varrho(n)$ of non-congruent regular integers modulo~$n$, corresponding to
sequence \seqnum{A055653} in the On-Line Encyclopedia of Integer Sequences
(OEIS), where an integer $m$ is said to be \emph{regular modulo~$n$} if the
congruence $m^2 x \equiv m \pmod{n}$ has a solution $x\in\mathbb{Z}$.  To
illustrate the significance of the sequence and Morgado's formula, we relate
them to a recent multi-prime, multi-power generalization of the RSA
cryptosystem.
\end{abstract}

\section{Introduction}

This work is motivated by a recent multi-prime, multi-power generalization of
the RSA cryptosystem \cite{Dohmen-Lange2025}, where the modulus is an
arbitrary integer $n>1$ and the messages are regular integers modulo~$n$. The
number of regular integers $m$ in $\mathbb{Z}_n = \{0,\dots,n-1\}$ is crucial
for estimating the probability of correct decryption in this generalized
scheme for random messages $m$ from the larger message space $\mathbb{Z}_n$.

The notion of a regular integer modulo $n$, defined below, dates back to
Morgado~\cite{Morgado1972}. Throughout, we let $\mathbb{N}$ denote the set of
positive integers.

\begin{definition}[\cite{Morgado1972}]
For each $n\in\mathbb{N}$, an integer $m$ is said to be \emph{regular
  modulo~$n$} if the congruence $m^2 x \equiv m \pmod{n}$ has a solution
$x\in\mathbb{Z}$.
\end{definition}

For each $n\in\mathbb{N}$, we use $\znr$ to denote the set of all
$m\in\mathbb{Z}_n$ that are regular modulo~$n$, and $\varrho(n)$ to denote its
cardinality.  The sequence $(\varrho(n))_{n\ge 1}$ appears as sequence
\seqnum{A055653} in the On-Line Encyclopedia of Integer Sequences (OEIS)
\cite{OEIS2023}. It was first studied by Morgado \cite{Morgado1972} and has
since been investigated by several
authors~\cite{Alkam2008,Apostol2013,Toth2008}.

Recall that $d\in\mathbb{N}$ is called a \emph{unitary divisor} of~$n$ if $d$
divides $n$ and $\gcd(d, n/d)=1$. Following Morgado \cite{Morgado1972}, we
write $d\udiv n$, if $d$ is a unitary divisor of~$n$.  Our focus is on the
following formula, due to Morgado \cite{Morgado1972}, in which $\varphi$
denotes Euler's totient function.

\begin{theorem}[\cite{Morgado1972}] For every $n\in\mathbb{N}$,
\begin{align}
  \label{eq:2}
  \varrho(n) & = \sum_{d\udiv n} \varphi(d) .
\end{align}
\end{theorem}

In this paper, we provide four proofs of this formula.  Unlike previously
published proofs \cite{Alkam2008,Toth2008}, the proofs presented here
are combinatorial in nature and do not rely on the multiplicativity
of~$\varrho$.  Instead, we repeatedly use the bijection principle and, in our
final proof, the inclusion-exclusion principle.  Continuing along this line of
reasoning, the multiplicativity of $\varrho$ follows naturally
from~\eqref{eq:2}.

The paper is organized as follows.  In Section~\ref{sec:preliminaries} we
provide a concise proof of Morgado's \cite{Morgado1972} characterization of
regular integers modulo~$n$, which we use in our proofs of Morgado's formula
\eqref{eq:2} in
Sections~\ref{sec:proof-equiv-relat}--\ref{sec:proof-incl-excl}. Each of our
four combinatorial proofs is self-contained and sheds a different light on the
formula.  From the authors' perspective, the purely bijective proof in
Section~\ref{sec:proof-bijection} is particularly noteworthy, as it yields an
encoding of the regular integers modulo~$n$ and may provide further insight
into the study of the sequence \seqnum{A055653}.

In Section~\ref{sec:application-rsa}, we relate this sequence and Morgado's
formula to the probability of correct decryption of a random message
$m\in\mathbb{Z}_n$ in a multi-prime, multi-power generalization of the RSA
cryptosystem, recently established by the present
authors~\cite{Dohmen-Lange2025}.
 
\section{Preliminaries}
\label{sec:preliminaries}

The following proposition, which is due to Morgado \cite{Morgado1972},
provides necessary and sufficient conditions for an integer $m$ to be regular
modulo $n$.  To keep this paper self-contained, we provide our own concise
proof.

\begin{proposition}[\cite{Morgado1972}]
\label{sec:regul-numb-unit-2}
For all $n\in\mathbb{N}$ and $m\in\mathbb{Z}$, the following statements
are equivalent:
\begin{enumerate}[(a)]
\item $m$ is regular modulo $n$,
\item $\gcd(m^2,n) = \gcd(m,n)$,
\item $\gcd(m,n)\udiv n$.
\end{enumerate}
\end{proposition}

\begin{proof}
$\text{(a)} \Leftrightarrow \text{(b)}$: In general, for every
$a,b\in\mathbb{Z}$, the congruence $ax\equiv b \pmod{n}$ has a solution
$x\in\mathbb{Z}$ if and only if $\gcd(a,n)\mid b$. Hence, $m$ is regular
modulo $n$ if and only if $\gcd(m^2,n)\mid m$, which in turn holds if and only
if $\gcd(m^2,n)=\gcd(m,n)$.
\par
$\text{(b)} \Rightarrow \text{(c)}$: Let $d=\gcd(m,n)$ and $g=\gcd(d,n/d)$.
By induction on $k$ we show that $g^k\mid d$ for all $k\ge 0$, which implies
$g=1$ and thus (c).  The case $k=0$ is trivial. For the induction step, assume
that $k>0$ and $g^{k-1}\mid d$. Then,
$g^k\mid dg\mid \gcd(m^2,n) = \gcd(m,n)=d$.
\par
$\text{(c)} \Rightarrow \text{(a)}$: The assumption states that
$\gcd(m,n/\gcd(m,n))=1$.  By Bezout's lemma, there exist $x,y\in\mathbb{Z}$
such that $m^2x + \frac{nm}{\gcd(m,n)}y = m$, whence $m^2x \equiv m \pmod{n}$.
\end{proof}

\section{Proof by equivalence relation}
\label{sec:proof-equiv-relat}

Our first proof of (\ref{eq:2}) is inspired by Morgado's original proof
\cite{Morgado1972}, but is considerably more formal and combinatorial, as it
makes explicit use of an equivalence relation and the bijection principle on
the resulting equivalence classes. Recall that
$\mathbb{Z}_n^\ast = \{m\in\mathbb{Z}_n\mid \gcd(m,n)=1\}$.

\begin{proof}
For $m_1,m_2\in\znr$, let $m_1\sim m_2$ if $\gcd(m_1,n)=\gcd(m_2,n)$; this
defines an equivalence relation on $\znr$.
By Proposition~\ref{sec:regul-numb-unit-2}, $m\in\znr$ if and only if
$\gcd(m,n)\udiv n$, so the equivalence classes are of the form $C_{n,d}$ with
$d\udiv n$, where
\begin{gather*}
C_{n,d} = \left\{ m\in\znr\mathrel| \gcd(m,n)=d \right\}.
\end{gather*}
Using the bijection principle, we show that for every unitary divisor $d$
of~$n$, 
\begin{gather}
\label{cnd}
|C_{n,d}| = |\mathbb{Z}_{n/d}^\ast| .
\end{gather}
To this end, define $h_{n,d} : C_{n,d}\rightarrow \mathbb{Z}_{n/d}^\ast$ by
$h_{n,d}(m) = m \bmod (n/d)$. This map is well defined,
since $\gcd(m \bmod (n/d),n/d) = \gcd(m,n/d) = \gcd(m,n,n/d) = \gcd(d,n/d)=1$.
It remains to show that $h_{n,d}$ is bijective.
\par
\emph{Injectivity.}  Suppose that $h_{n,d}(m_1) = h_{n,d}(m_2)$. Then,
$m_1 \equiv m_2 \pmod{n/d}$.  Since $m_1,m_2\in C_{n,d}$, we have
$\gcd(m_1,n) = \gcd(m_2,n)=d$, which implies $m_1 \equiv m_2 \pmod{d}$.
Because $d$ and $n/d$ are coprime, combining both congruences gives
$m_1\equiv m_2 \pmod{n}$, and hence $m_1=m_2$.
\par
\emph{Surjectivity.} Let $d'\in\mathbb{Z}_{n/d}^\ast$, and define
$m = d ((d'\/i) \bmod (n/d))$, where $i$ denotes an inverse of $d$ modulo
$n/d$. We claim:
\begin{enumerate}[(i)]
\item \label{first} $m\in C_{n,d}$\,;
\item \label{second} $h_{n,d}(m)=d^\prime$.
\end{enumerate}
For \eqref{first}, it suffices to show that $\gcd(m,n)=d$. Indeed, since
$d\mid m$ and $d\mid n$, and since $\gcd(d',n/d)=1$ and $\gcd(i,n/d)=1$, we have
\[ \gcd(m,n) = d \gcd(d'\/i \bmod (n/d), n/d) = d \gcd(d'\/i, n/d) = d \gcd(i,
n/d) = d. \]
Part~\eqref{second} follows immediately, since
\[ h_{n,d}(m) = m \bmod (n/d) = (d\/i  \bmod (n/d))(d' \bmod (n/d))=d'. \]
From \eqref{cnd} and the disjointness of the equivalence
classes, we conclude that
\begin{gather*}
\label{eq:8}
\varrho(n) = \sum_{d\udiv n} \, |C_{n,d}| = \sum_{d\udiv n} \, |\mathbb{Z}_{n/d}^\ast| =\sum_{d\udiv n} \,
\varphi\left( \frac{n}{d} \right) = \sum_{d\udiv n} \varphi(d) ,
\end{gather*}
which proves~\eqref{eq:2}.
\end{proof}

\section{A purely bijective proof}
\label{sec:proof-bijection}

Our next proof is purely bijective.  The idea is to establish a bijection
between $\znr$ and the set of pairs $(d,d')$ with $d\udiv n$ and
$d'\in\mathbb{Z}_d^\ast$ that are counted by the right-hand side of
(\ref{eq:2}). This bijection yields an encoding of $\znr$ that may prove useful
beyond this proof.

\begin{proof}
Let $U_n$ denote the set of unitary divisors of~$n$.
Consider the map
\begin{gather*}
\label{eq:3}
f_n : \znr \rightarrow \left\{ (d,d')\mathrel| d\in U_n,
d'\in\mathbb{Z}_d^\ast \right\},
\end{gather*}
defined by
\begin{gather*}
\label{eq:6}
f_n(m) := \left( \frac{n}{\gcd(m,n)}, m\bmod \frac{n}{\gcd(m,n)} \right) .
\end{gather*}
We first show that $f_n$ is well defined.  Let $d = n/\gcd(m,n)$.  Then
$d\in U_n$, and hence
$\gcd(d, m \bmod d) = \gcd(d, m) = \gcd(d,m,n) = \gcd(d,\gcd(m,n)) =
\gcd(d,n/d)=1$, which implies $m\bmod d\in\mathbb{Z}_d^\ast$. To apply
the bijection principle, we show that $f_n$ is bijective.
\par
\emph{Injectivity.} Suppose that $f_n(m_1) = f_n(m_2)$. Then,
$\gcd(m_1,n) = \gcd(m_2,n)$, which we denote by~$d$.  Evidently,
$m_1\bmod n/d = m_2\bmod n/d$, which means that $m_1\equiv m_2 \pmod{n/d}$.
From $\gcd(m_1,n)=d$ we can write $m_1=d m_1'$, $m_2=d m_2'$, and $n=d n'$
with $\gcd(m_1',n')=\gcd(m_2',n')=1$. Therefore, $m_1-m_2 = d(m_1'-m_2')$, so
$d\mid m_1-m_2$, which gives $m_1\equiv m_2 \pmod{d}$.  Since $d$ and $n/d$
are coprime (because $d\udiv n$), combining the congruences
$m_1\equiv m_2\pmod{d}$ and $m_1\equiv m_2\pmod{n/d}$ gives
$m_1\equiv m_2 \pmod n$, and hence $m_1=m_2$.
\par
\emph{Surjectivity.} Let $d\in U_n$ and $d'\in\mathbb{Z}_d^\ast$.  We define
$m$ as
\begin{gather}
\label{eq:1}
m = \frac{n}{d} \big(( d'\/j ) \bmod d \big),
\end{gather}
where $j$ is an inverse of $n/d$ modulo~$d$. We claim:
\begin{enumerate}[(i)]
\item \label{firstclaim} $m\in\znr$\,;
\item \label{secondclaim} $f_n(m)=(d,d')$.
\end{enumerate}
Since $n/d$ divides both $m$ and $n$, and since $\gcd(d',d)=1$ and
$\gcd(j,d)=1$, we have
\begin{gather}
\label{eq:5}
\gcd(m,n)  = \frac{n}{d} \gcd\left( \left( d'\/j \right) \bmod d,
            d\right) = \frac{n}{d} \gcd\left( d'\/j ,
            d\right) = \frac{n}{d} \gcd\left( j,
            d\right) = \frac{n}{d}.
\end{gather}
Hence $\gcd(m,n)\udiv n$, and by Proposition~\ref{sec:regul-numb-unit-2},
$m\in\znr$, as claimed in \eqref{firstclaim}.  For part \eqref{secondclaim}, we
note that $d = n/\gcd(m,n)$ follows from \eqref{eq:5}, and $d'=m\bmod d$
follows from \eqref{eq:1}, since $j$ is an inverse of $n/d$ modulo~$d$.  Thus,
\eqref{firstclaim} and \eqref{secondclaim} are shown, and the proof is
complete.
\end{proof}

To illustrate the proof, we list the assignments $m\mapsto f_{20}(m)$ for
$m\in\mathbb{Z}_{20}^{\text{reg}}$:
\begin{align*}
0 & \mapsto (1,0), & 4 & \mapsto(5,4), & 8 & \mapsto (5,3), & 12 & \mapsto (5,2), &
                                                             16 & \mapsto
                                                                 (5,1), \\
1& \mapsto (20,1), & 5 & \mapsto (4,1), & 9 & \mapsto (20,9), & 13 & \mapsto
                                                                     (20,13),
                                                                                   &
                                                                                     17 & \mapsto (20,17),\\  
3 & \mapsto (20,3), & 7 & \mapsto (20,7), & 11 & \mapsto (20,11),
& 15 & \mapsto (4,3), & 19 & \mapsto (20,19).
\end{align*}  

\begin{remark}
In view of \eqref{eq:1}, the inverse of $f_n$ takes the form $f_n^{-1}(d,d') =
\frac{n}{d} ( ( ( n/d \bmod d)^{-1} d' ) \bmod d )$ for every $d\in U_n$ and
$d'\in\mathbb{Z}_d^\ast$.
\end{remark}

\begin{remark}
The proof can be restated by defining
$f_n(m) := ( \gcd(m,n), m\bmod n/\gcd(m,n))$, which maps
from $\znr$ to $\{ (d,d')\mathrel| d\in U_n, d'\in\mathbb{Z}_{n/d}^\ast \}$.
In this setting, $f_n^{-1}(d,d') = d(((d \bmod (n/d))^{-1}d')\bmod (n/d))$.
\end{remark}

\section{Proof by reduced fractions}
\label{sec:proof-reduced-fractions}

Our third proof is inspired Gauss's formula $n = \sum_{d\mid n} \varphi(d)$,
as reproduced in the textbook by Graham, Knuth, and Patashnik
\cite[pp.~134--135]{Graham1994}.  The key idea is to establish a bijection
between $\znr$ and the set of reduced fractions of the form $k/d$, where
$d\udiv n$ and $k<d$.

\begin{proof}
Consider the fractions $m/n$ with $m\in\znr$. Reducing these fractions to
lowest terms yields fractions of the form
\[ k/d = (m/\gcd(m,n))/(n/\gcd(m,n)). \] By
Proposition~\ref{sec:regul-numb-unit-2}, $m$ is regular modulo $n$ if and only
if $\gcd(m,n)\udiv n$, or equivalently, if and only if $n/\gcd(m,n)\udiv
n$. Hence, the denominators of these reduced fractions are precisely the
unitary divisors $d$ of $n$. Each reduced fraction $k/d$ with $d\udiv n$ and
$k<d$ arises in this way by reducing $(kn/d)/n$ to lowest terms. To
complete the argument, we show that $kn/d\in\znr$.  Because $k$ and $d$, as well
as $d$ and $n/d$, are coprime, 
\[ \gcd(kn/d,n) = \gcd(k(n/d),d(n/d)) = n/d \udiv n. \]
Hence, by
Proposition~\ref{sec:regul-numb-unit-2}, $kn/d\in\znr$.  Thus, the
$\varrho(n)$ reduced fractions can be grouped according to their denominator
$d\udiv n$, with $\varphi(d)$ reduced fractions for each denominator~$d$.
\end{proof}

To illustrate the proof, consider the $\varrho(20)$ fractions $m/20$ for
$m\in\mathbb{Z}_{20}^{\text{reg}}$:
\begin{align*}
\frac{0}{20} &\,,&\!\! \frac{1}{20} &\,,&\!\! \frac{3}{20} &\,,&\!\! \frac{4}{20} &\,,&\!\! \frac{5}{20}
  &\,,&\!\! \frac{7}{20} &\,,&\!\! \frac{8}{20} &\,,&\!\! \frac{9}{20} &\,,&\!\! \frac{11}{20} &\,,&\!\! \frac{12}{20} &\,,&\!\!
                                                                      \frac{13}{20}
  &\,,&\!\! \frac{15}{20} &\,,&\!\! \frac{16}{20} &\,,&\!\! \frac{17}{20} &\,,&\!\!
                                                                    \frac{19}{20} \, . \\
  \intertext{Grouping the reduced fractions by their denominators yields}
  \frac{0}{1} &\,,&\!\! \frac{1}{4} &\,,&\!\! \frac{3}{4} &\,,&\!\! \frac{1}{5} &\,,&\!\! \frac{2}{5}
  &\,,&\!\! \frac{3}{5} &\,,&\!\! \frac{4}{5} &\,,&\!\! \frac{1}{20} &\,,&\!\! \frac{3}{20} &\,,&\!\!
                                                                        \frac{7}{20}
                                                                                        &\,,&\!\!
                                                                                             \frac{9}{20}
  &\,,&\!\! \frac{11}{20} &\,,&\!\! \frac{13}{20} &\,,&\!\! \frac{17}{20} &\,,&\!\!
                                                                    \frac{19}{20} \, ,
\end{align*}
with $\varphi(1)=1$ fraction having denominator $1$,
$\varphi(4)=2$ fractions having denominator $4$,
$\varphi(5)=4$ fractions having denominator $5$,
and $\varphi(20)=8$ fractions having denominator~$20$.
Hence $\varrho(20) = 1 + 2 + 4 + 8 = 15$.

\begin{remark}
There is an obvious connection with the proof in
Section~\ref{sec:proof-bijection}: a fraction $a/b$ appears in the list of
reduced fractions if and only if $f_n(m) = (b,a)$ for some $m\in\znr$.
\end{remark}

\section{Proof by inclusion-exclusion}
\label{sec:proof-incl-excl}

Our final proof of \eqref{eq:2} is based on a combined application of the
inclusion-exclusion principle, the bijection principle, and the
multiplicativity of Euler's totient function $\varphi(n)$.

\begin{proof}
For every integer $m\ge 0$ and every prime $p$, let $\nu_p(m)$ denote the
multiplicity of $p$ in the prime factorization of~$m$.  For every
$m\in\mathbb{Z}_n$, we have $m\in\znr$ if and only if $\nu_p(m)=0$ or
$\nu_p(m)\ge \nu_p(n)$ for each prime divisor $p$ of $n$, as follows from
Proposition~\ref{sec:regul-numb-unit-2}. Let $P(n)$ denote the set of prime
divisors of $n$, and for each $p\in P(n)$, define
\begin{gather*}
A_p  =  \left\{m\in\mathbb{Z}_n\mathrel| 0<\nu_p(m)<\nu_p(n) \right\} .
\end{gather*}
Then by the inclusion-exclusion principle,
\begin{gather}
\label{eq:ie1}
 \varrho(n) = \left    | \bigcap_{p\in P(n)} \overline{A_p} \right| = \sum_{I\subseteq P(n)} (-1)^{|I|} \left| \bigcap_{i\in I} A_i \right| .
\end{gather}
In this formula, $m\in\bigcap_{i\in I} A_i$ if and only if
$m = k \prod_{i\in I} i$ for some $k\le \frac{n}{\prod_{i\in I}i}$ such that
$j^{\nu_j(n)-1} \nmid k$ for each $j\in I$; that is, if and only if
$k\in\bigcap_{j\in I} \overline{B_j}$, where
\[ B_j = \left\{ \left. 1\le k\le \frac{n}{\prod_{i\in I} i} \,\,\right|\,
j^{\nu_j(n)-1} \mid k \right\} \quad (j\in I). \] Clearly,
$m \mapsto \frac{m}{\prod_{i\in I} i}$ defines a bijection from
$\bigcap_{i\in I} A_i$ to $\bigcap_{j\in I} \overline{B_j}$.  Therefore, by
the bijection principle and another application of the inclusion-exclusion
principle, we have
\begin{gather}
\label{eq:ie3}
\left| \bigcap_{i\in I} A_i \right|
= \sum_{J\subseteq I} (-1)^{|J|} \left| \bigcap_{j\in J} B_j \right| = \sum_{J\subseteq I} (-1)^{|J|} \frac{n}{\prod_{i\in I}i \prod_{j\in J}
  j^{\nu_j(n)-1}}  .
\end{gather}
Combining (\ref{eq:ie1}) and (\ref{eq:ie3}) and then changing the order of
summation, we obtain
\begin{gather*}
\varrho(n) = \sum_{I\subseteq P(n)} \sum_{J\subseteq I} (-1)^{|I|+|J|}
             \frac{n}{\prod_{i\in I\setminus J} i \prod_{j\in J} j^{\nu_j(n)}}
  = \sum_{J\subseteq P(n)} \prod_{j\in J} \frac{n}{j^{\nu_j(n)}}
  \sum_{I\supseteq J} (-1)^{|I|+|J|} \prod_{i\in I\setminus J} \frac{1}{i} .
\end{gather*}
Replacing $J$ by its complement in $P(n)$, and factoring the
inner sum, it follows that
\begin{gather*}
\varrho(n) = \sum_{J\subseteq P(n)} \prod_{j\in J}  j^{\nu_j(n)}   \sum_{I\subseteq
  J} (-1)^{|I|} \prod_{i\in I} \frac{1}{i} \notag  = \sum_{J\subseteq P(n)} \prod_{j\in J}
j^{\nu_j(n)}
\left( 1 - \frac{1}{j} \right).
\end{gather*}
Using Euler's totient function and its multiplicativity, we obtain
\begin{gather*}
\varrho(n) = \sum_{J\subseteq P(n)}
\prod_{j\in J} \varphi(j^{\nu_j(n)}) = \sum_{J\subseteq P(n)} \varphi\left( \prod_{j\in J}
j^{\nu_j(n)} \right) .
\end{gather*}
We finally observe that the last sum ranges over all positive divisors
$d = \prod_{j\in J} j^{\nu_j(n)}$ of $n$ that are coprime to $n/d$, i.e., over
all unitary divisors $d$ of $n$, thus proving \eqref{eq:2}.
\end{proof}

\section{Significance for cryptography}
\label{sec:application-rsa}

The authors~\cite{Dohmen-Lange2025} encountered regular integers modulo $n$
while developing a generalization of the RSA scheme \cite{RSA1978} to
arbitrary multi-prime, multi-power moduli.  For such a generalized modulus
$n=p_1^{e_1}\dots p_r^{e_r}$ with distinct primes $p_1,\dots,p_r$ and
exponents \mbox{$e_1,\dots,e_r\in\mathbb{N}$}, the public key $(n,e)$ and the
private key $(n,d)$ are established in the same way as in the classical RSA
scheme: choose $1<e<\varphi(n)$ such that $\gcd(e,\varphi(n))=1$, and compute
$1<d<\varphi(n)$ such that $ed\equiv 1\pmod{\varphi(n)}$.  As in classical
RSA, a message $m\in\mathbb{Z}_n$ is encrypted by raising $m$ to the $e$-th
power modulo~$n$ and decrypted by raising $m$ to the $d$-th power modulo~$n$.

A key observation, proved by the present authors \cite{Dohmen-Lange2025}, is
that decryption reverses encryption if and only if the message is regular
modulo~$n$. Consequently, by \eqref{eq:2}, the probability of correct
decryption of a random message from $\mathbb{Z}_n$ is given by
\begin{align*}
  \frac{\varrho(n)}{n} & = \frac{1}{n} \sum_{d\udiv n} \varphi(d) ,
\end{align*}                         
which illustrates the significance of the sequence \seqnum{A055653} and
Morgado's formula \eqref{eq:2} in the context of cryptography.  As further
shown by the present authors~\cite{Dohmen-Lange2025},
\begin{align*}
  \frac{\varrho(n)}{n} & \ge 1-\frac{r}{2^{k-1}} \, ,
\end{align*}
where $n = p_1^{e_1}\dots p_r^{e_r}$ with distinct $k$-bit primes
$p_1,\dots,p_r$.  Therefore, even for today's standard choices of~$k$, for
example $k=1024$, almost all messages in $\mathbb{Z}_n$ are decrypted
correctly, and the restriction to regular messages is negligible.  Although this
conclusion is satisfactory from a practical point of view, there remains
potential for sharper bounds on the correctness probability.  Asymptotic
results on $\varrho(n)$ and related quantities such as
$\varrho(n)/\varphi(n)$, as obtained by {Apostol} and {Petrescu}
\cite{Apostol2013} and by {T\'oth} \cite{Toth2008}, may prove crucial in this
regard.

\end{document}